\documentclass[11pt]{article} \usepackage{amsfonts}   

\makeatletter
\def\USquality{\textheight=250mm \textwidth=170mm \topmargin=0Truein
               \voffset=-1truein
               \ifcase \@ptsize \hoffset=-23mm
                       \or \hoffset=-20mm \or \hoffset=-15mm \fi}
\makeatother
\USquality 

\def\beq#1#2{\begin{equation} \label{#1} #2 \end{equation}}
\def\bea#1{\begin{eqnarray*} #1 \end{eqnarray*}} \def\a{\!\!\!&\!\!\!\!&}

\def\thname{Theorem}     \def\lmname{Lemma}      \def\prname{Proposition}
\def\dfname{Definition}  \def\crname{Corollary}  \def\rmname{Remark}

\newtheorem{theorem}{\thname}[section]   
\newtheorem{lemma}{\lmname}[section]     
\newtheorem{corollary}[lemma]{\crname}   

\newtheorem{dftn}{\dfname}[section]
\newenvironment{definition}{\begin{dftn}\rm}{\end{dftn}} 
\newtheorem{rmrk}[lemma]{\rmname}

             \catcode`@=11 \@addtoreset{equation}{section} \catcode`@=12

\def\proof{\smallskip \noindent {\bf Proof. \ }}       
\def\blanksquare{\,\,\,$\sqcup\!\!\!\!\sqcap$}         
\def\qed{\hfill\blanksquare\linebreak\smallskip\par}   

\def\function#1{\left\{\!\!\!\begin{array}{ll} #1 \end{array} \right.}

\newcommand\mlbscale{1pt} 
\newif\iffigs\figstrue 

\def\bfig(#1,#2)#3#4{\begin{figure} \begin{center}
    \framebox{\setlength{\unitlength}{\mlbscale} 
       \iffigs \begin{picture}(#1,#2) #3 \end{picture} 
       \else \begin{picture}(60,10)(0,0) 
                   \put(0,0){\framebox(60,10){Figure}} \end{picture} \fi} 
    \end{center} \caption{#4} \end{figure}}

\def\Bfig(#1,#2)#3#4{\begin{figure} \begin{center}
    \setlength{\unitlength}{\mlbscale} 
       \iffigs \begin{picture}(#1,#2) #3 \end{picture} 
       \else \begin{picture}(60,10)(0,0) 
                   \put(0,0){\framebox(60,10){Figure}} \end{picture} \fi 
    \end{center} \caption{#4} \end{figure}}

\def\bpic(#1,#2)#3{\setlength{\unitlength}{\mlbscale} 
    \begin{picture}(#1,#2) #3 \end{picture}} 

\def\bline(#1,#2)(#3,#4)(#5){\put(#1,#2){\line(#3,#4){#5}}}  

\def\ep{\varepsilon} \def\la{\lambda} \def\phi{\varphi}

\def\map{T} \def\n{\noindent} \def\supp{{\rm supp}}
\def\IR{{\mathbb{R}}} \def\IZ{{\mathbb{Z}}}
\def\cB{{\cal B}} \def\cN{{\cal N}} \def\cmap{{\vec{\cal T}}}
\def\cM{{\cal M}} \def\cI{{\vec{\cal I}}} \def\La{\Lambda}
         \def\*#1{#1^*}
\def\toas#1{\stackrel{#1}{\longrightarrow}}

\def\OPR#1{\begin{definition} #1 \end{definition}}

\begin{document}
\title{Long range action in networks of chaotic elements}
\author{Michael Blank\thanks{Russian Academy of Sciences,
                             Inst. for Information Transmission Problems,
                             and Observatoire de la Cote d'Azur, ~
                             e-mail: blank@obs-nice.fr}, ~
        Leonid Bunimovich\thanks{School of Mathematics,
                                 Georgia Institute of Technology;
                                 e-mail: bunimovh@math.gatech.edu}}
\date{\today}
\maketitle

{\bf Abstract}. We show that under certain simple assumptions on
the topology (structure) of networks of strongly interacting chaotic
elements a phenomenon of long range action takes place, namely
that the asymptotic (as time goes to infinity) dynamics of an
arbitrary large network is completely determined by its boundary
conditions. This phenomenon takes place under very mild and
robust assumptions on local dynamics with short range interactions.
However, we show that it is unstable with respect to arbitrarily
weak local random perturbations.

\section{Introduction}\label{s:intro}
Despite a number of impressive results in the field of 
coupled map lattices (CML) up to now only the situation when the
interaction between individual systems is very weak is reasonably
well understood from the mathematical point of view (see a review
of recent results, e.g. in \cite{Bu,KL}). In that case it is
known that under mild technical assumptions the asymptotic
behavior of a CML is very close to the behavior of a direct
product of individual uncoupled systems. Note that
counterexamples demonstrating very different types of behavior
under arbitrary weak interactions are known as well (see \cite{BB}
for details). On the other hand, numerical studies (see, e.g.
\cite{MH}) suggest that when the interactions becomes stronger a
kind of phase transition from one to several invariant SBR
measures takes place. Despite a number of efforts during last ten
years to give a mathematical description of this phenomenon no
rigorous results for infinite systems have been obtained.
Probably one of the reasons is that all mathematical results in
the field are of perturbative nature and corresponding methods do 
not work well for large perturbations. Another problem here is the 
necessity to work with infinite collections of individual maps. 
Therefore, it is harder to study corresponding functional spaces.

In this paper we shall follow a more ``physical'' approach to the
analysis of CMLs' dynamics in a strong interactions case. Indeed,
a general receipt in statistical physics to an analysis of an
infinitely large system is to consider its restriction to a
finite volume (box) with certain boundary conditions and to study
what happens when the size of the box grows to infinity. The case
when the dynamics inside of the (sufficiently large) box is not
influenced by boundary conditions corresponds to the situation known
as the case (range) of weak interactions (or as the absence of
phase transitions). We shall be mainly interested in the
opposite situation when the asymptotic (as time goes to infinity)
dynamics inside of the box (with an arbitrary large diameter)
is completely determined by the boundary conditions and which is
reasonably to call the long range action (LRA). As we shall
see the LRA takes place under sufficiently strong interactions in 
networks whose topologies satisfy the property which we call
{\em unidirectional} and which generalizes well known
one-dimensional CMLs with unidirectional interactions \cite{PRK}.

The paper is organized as follows. In the next section we give 
some technical definitions and formulate our main result about
the existence of the LRA in a reasonably large class of CMLs
(Theorem~\ref{t:LRA}). Sections~\ref{s:LRA-proof} and \ref{s:enum}
contain the proof of this result and the proof of the equivalence
of the notion of the unidirectional interaction to the absence of
cycles in the directed connectivity graph of interactions in
a network. Section~\ref{s:bc} is dedicated to the discussion of
the topology of networks corresponding to CMLs under study and
to the influence of different types of boundary conditions on the
existence of LRA. In Section~\ref{s:random} we study random
perturbations of local dynamics and show that the LRA phenomenon
is stochastically unstable. Finally in Section~\ref{s:conclusion}
we briefly address some open questions and compare our results to
some situations known in the theory of synchronization.

\section{The setup}\label{s:setup}

Let $\map:X\to X$ be a map from a compact metric space $(X,\rho)$
into itself satisfying the Lipschitz condition with the contracting
constant $\La_\map$, i.e. %
\beq{e:contr-const}{\Lambda_{\map} := \sup_{x\ne y\in X}
                    \frac{\rho(\map x, \map y)}{\rho(x,y)}<\infty,}%
(e.g. an Anosov map on a finite-dimensional torus with the Euclidean
metric). Then the pair $(\map,X)$ defines a dynamical system.
Consider an infinite collection of dynamical systems $(\map_i,X)$
with $i\in\IZ^d$ and a certain interaction between them. We always
assume that
$$ \Lambda_{\map_i} \le \Lambda_{\map}<\infty $$
for all $i\in\IZ^d$ and denote the direct product of the systems
$(\map_i,X)$ by $(\vec{\map},\vec{X})$. The maps $\map_i$ we shall
call {\em local} maps. Elements $\vec{x}$ of the infinite
dimensional space $\vec{X}$ will be called {\em vectors} and
their coordinates will be denoted as $x_i:=(\vec{x})_i\in X$ with
$i\in\IZ^d$.

\OPR{By the {\em interaction} on the lattice
$\vec{X}:=X^{\IZ^d}$ we
shall mean the map $\cI:\vec{X}\to\vec{X}$ preserving the
`diagonal',  i.e. $\cI(C\vec{e})=C\vec{e}$ for any constant $C$
and the vector $\vec{e}\in\vec{X}$ with all entries equal to $1$.
We restrict ourselves to the case of continuous interactions
satisfying the assumption that there exists a constant $\La_\cI<1$
such that%
\beq{e:lip-int}{\rho((\cI\vec{x})_\ell, (\cI\vec{y})_\ell)
                \le \La_\cI \sum_{i\in\IZ^d} \rho(x_i,y_i) } %
for any $\ell\in\IZ^d$ and any vectors $\vec{x},\vec{y}\in\vec{X}$.}

In particular, this means that for any pair of vectors which differ
at only one coordinate the map $\cI$ is a pure contraction.
It is tempting to assume that the interaction $\cI$ is a pure
contraction in the product metric, which would simplify a lot the
further analysis. At the first sight it seems that the difference
between the assumption~(\ref{e:lip-int}) and the pure contraction
is not essential. However a closer look shows that the latter
contradicts to the definition of the interaction which says that
the map $\cI$ preserves the `diagonal', i.e. the distance between
the vectors $\cI(C_1\vec{e})$ and $\cI(C_2\vec{e})$ is the same as
the distance between the vectors $C_1\vec{e}$ and $C_2\vec{e}$ for
any constants $C_1,C_2$.

The {\em strength} of the interactions we shall measure by the
difference $1-\Lambda_{\cI}$. Indeed, without the interactions one
has $\Lambda_{\cI}=1$. A ``typical'' example: $\cI$ is a linear
map represented by an infinite matrix with normalized rows. Local
interactions correspond to the case when only terms near the diagonal
of the matrix differ from zero. Again it is useful to have in mind
a ``standard'' example of the local interaction (often called a
``diffusive coupling'') on a one-dimensional lattice when after the
interactions the new value at the site $i$ is calculated according to
the formula: %
\beq{e:diff-coupling}{ x_i \to (1-2c)x_i + c(x_{i-1} + x_{i+1}) .}%
Here the interaction is parametrized by the value $0\le c\le1/2$
and $\Lambda_\cI=\max(1-2c,c)$.

Another important and well known example describes unidirectional
interactions on a one-dimensional lattice defined by the formula: %
\beq{e:lin-unidir}{ x_i \to (1-c)x_i + cx_{i-1} .} %
In this case $0\le c\le1$ and $\Lambda_\cI=\max(1-c,c)$.

Both these examples belong to the class of {\em short range} 
interactions when subsystems interact only if the Eucledian distance 
between the corresponding sites do not exceed a certain threshold.
To emphasize that our results do not come as a consequence of long 
range interactions (when subsystems being arbitrary far apart may 
interact with each other) we shall restrict the consideration only 
to short range interactions.

\OPR{The superposition $\cmap:=\cI\circ\vec{\map}$
of the interaction and the direct product of local maps defines
the dynamical system $(\cmap,\vec{X})$ which is called
a {\em coupled map lattice} (CML).}

\OPR{A {\em box} $\cB$ is a bounded subset of the integer
lattice $\IZ^d$.}

To restrict the dynamics to the box $\cB$ we define it as
follows: %
\beq{e:map-in-box}{
  \left(\cmap_{|\cB}~\vec{x}\right)_i :=
   \function{
     \left(\left(\cI\circ\vec{\map}\right)\vec{x}\right)_i
           &\mbox{if } i \in \cB\\
     x_i    &\mbox{otherwise}.} }%

\OPR{A collection of values of $\vec{x}$ outside of the box $\cB$
is called {\em boundary conditions} and is denoted by
$\vec{x}_{|\cB^c}$, where by $\cB^c$ we denote the complement to
the set $\cB$.}

If we choose certain boundary conditions the map $\cmap_{|\cB}$
defines a finite dimensional dynamics inside a finite box $\cB$.
The collection of functions
$x_i(t):=\left(\cmap_{|\cB}^t~\vec{x}(0)\right)_i$ with $i\in\cB$
and $t\in\IZ_+$ will be called a {\em solution} of our problem
corresponding to the initial data $\vec{x}(0)$. Clearly $\vec{x}(0)$
can be decomposed into the genuine initial data $\vec{x}_{|\cB}(0)$
and the boundary conditions $\vec{x}_{|\cB^c}$.

\OPR{By the {\em long range action} (LRA) we shall mean that for
a given box $\cB$ and given boundary conditions $\vec{x}_{|\cB^c}$
there exists a vector $\vec{\hat{x}}_{|\cB}$ (called a
{\em limit solution}) such that for any initial data
$\vec{x}_{|\cB}(0)$ the corresponding solution $\vec{x}(t)$
converges pointwise with time $t$ to $\vec{\hat{x}}_{|\cB}$. }

The main issue here is that the boundary conditions define the long
term dynamics inside the box $\cB$ in a unique way contrary to the
usual situation when the dynamics inside the box is asymptotically
independent on the boundary conditions. On the other hand, we do not
assume that different boundary conditions lead to different solutions.
The point is that if the local maps are not bijective a pair of very
different boundary conditions may lead to the same limit solution.
On the other hand, under some simple additional assumptions (of
non-degeneracy type) close but different boundary conditions lead
to different limit solutions.\footnote{See a discussion in the end
    of Section~\ref{s:LRA-proof}.}

\OPR{We shall say that the interaction $\cI$ in the box $\cB$ is
{\em unidirectional} if there is a positive function
$\cN:\cB\to\IZ^1$ such that \\
-- $\forall i\ne j\in\cB$ we have $\cN(i)\ne\cN(j)$;\\
-- $\forall \vec{x}\in\vec{X}$ and $\forall i\in\cB$ the
   value $(\cI\vec{x})_i$ does not depend on the values $x_j$ for
   all $j$ with $\cN(j)>\cN(i)$.}

In other words the function $\cN$ is an enumeration according to
which sites in $\cB$ with larger numbers do not `interact' with
sites with smaller numbers.

It is convenient to describe a topology of interactions in terms
of connectivity graphs.

\OPR{A {\em connectivity graph} $G:=G(\cB,\cI)$ of the
interaction $\cI$ in the box $\cB$ is a directed graph $G$, whose
vertices correspond to the sites of $\cB$. Besides the edge $(i,j)$
connecting the vertex $i\in G$ to the vertex $j\in G$ belongs to $G$
if and only if the value $(\cI\vec{x})_j$ depends on the value $x_i$.
The latter means that the value $(\cI\vec{x})_j$ assumes at least
two different values when we are changing $x_i$ assuming that all
other coordinates of $\vec{x}$ are fixed.}

To take into account the boundary conditions we include to the graph
$G$ all edges leading from $\cB^c$ to the points belonging to the
box $\cB$ and extend the function $\cN$ to these sites setting its
value to $0$ on them.

For each site $\ell\in\cB$ we denote by $\Gamma(\ell)$
the collection of sites in the box $\cB$ from which there are
directed paths in the graph $G$ leading to $\ell$, and denote by
$L(\ell)$ the length of the shortest path in $G$ leading to $\ell$
from the points outside the box $\cB$. If there are no such paths
we set $L(\ell)=\infty$.

Now we are ready to formulate the main result of the paper.

\begin{theorem}\label{t:LRA} Let the interaction in a finite box
$\cB$ be unidirectional and let for any $i\in\cB$ there be a path
in $G(\cB,\cI)$ from a site outside of the box $\cB$ to $i$. Then
for any boundary conditions the LRA takes place whenever
$\Lambda_\cI\Lambda_\map<1$. \end{theorem}

Note the local maps may be both regular and chaotic here. 
The simplest example of the unidirectional interaction is a linear
one-dimensional (i.e. $d=1$) system of unidirectionally interacting
maps with the interaction described by the
relation~(\ref{e:lin-unidir}). Of course, if this situation
would exhaust the list of examples there would be no reason to
discuss them at all. It turned out that there is a much bigger
class of systems with effectively unidirectional interactions.
This class of interactions have a very simple and intuitive
interpretation in terms of a connectivity graph.

\begin{theorem}\label{t:equiv} An interaction is unidirectional
if and only if the corresponding directed connectivity graph has
no cycles.
\end{theorem}

Two examples of graphs satisfying this property are depicted in
Fig.~\ref{f:graphs}. The left panel in Fig.~\ref{f:graphs}
corresponds to the interaction similar to the one which is
well known under the name Toom's North-East voting model
in statistical physics (see e.g. \cite{Li}), while the right one
demonstrates a bit more complex topology of interactions.
Observe that in both examples the graphs may be continued
periodically to an arbitrary large rectangular box.

\Bfig(320,150)
      {
       \bline(0,0)(1,0)(150)   \bline(0,0)(0,1)(150)
       \bline(0,150)(1,0)(150) \bline(150,0)(0,1)(150)
       \thicklines
       \bline(20,30)(1,0)(110)
       \bline(20,75)(1,0)(110)
       \bline(20,120)(1,0)(110)
       \bline(30,20)(0,1)(110)
       \bline(75,20)(0,1)(110)
       \bline(120,20)(0,1)(110)
       \put(100,120){\vector(-1,0){7}}
       \put(100, 75){\vector(-1,0){7}}
       \put(100, 30){\vector(-1,0){7}}
       \put( 60,120){\vector(-1,0){7}}
       \put( 60, 75){\vector(-1,0){7}}
       \put( 60, 30){\vector(-1,0){7}}
       \put( 30, 60){\vector(0,-1){7}}
       \put( 30,100){\vector(0,-1){7}}
       \put( 75, 60){\vector(0,-1){7}}
       \put( 75,100){\vector(0,-1){7}}
       \put(120, 60){\vector(0,-1){7}}
       \put(120,100){\vector(0,-1){7}}
       \put(160,0){\bpic(150,150){
         \bline(0,0)(1,0)(150)   \bline(0,0)(0,1)(150)
         \bline(0,150)(1,0)(150) \bline(150,0)(0,1)(150)
         \bline(20,30)(1,0)(110) \bline(20,120)(1,0)(110)
         \bline(30,20)(0,1)(110) \bline(120,20)(0,1)(110)
         \bline(30,90)(1,1)(30)  \put(30,90){\vector(1,1){15}}
         \bline(30,60)(1,-1)(30) \put(30,60){\vector(1,-1){15}}
         \bline(120,90)(-1,1)(30) \put(120,90){\vector(-1,1){15}}
         \bline(120,60)(-1,-1)(30) \put(120,60){\vector(-1,-1){15}}
         \put(30,90){\vector(0,1){20}}
         \put(30,60){\vector(0,-1){20}} \put(30,60){\vector(0,1){20}}
         \put(120,90){\vector(0,1){20}}
         \put(120,60){\vector(0,-1){20}} \put(120,60){\vector(0,1){20}}
         \put(60,30){\vector(-1,0){20}}  \put(60,30){\vector(1,0){20}}
                                         \put(90,30){\vector(1,0){20}}
         \put(60,120){\vector(-1,0){20}}  \put(60,120){\vector(1,0){20}}
                                          \put(90,120){\vector(1,0){20}}
                  }}
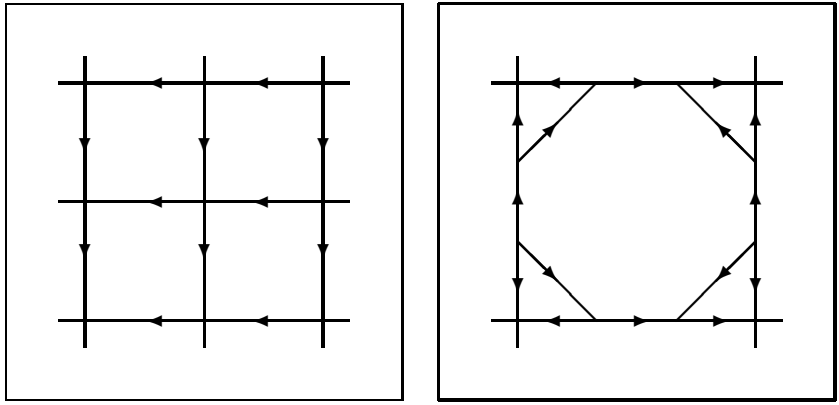
      }{Examples of unidirectional connectivity graphs \label{f:graphs}}

\section{Proof of Theorem~\ref{t:LRA}} \label{s:LRA-proof}

First we prove the following technical result.

\begin{lemma}\label{l:contraction} Let $\{v_i(t)\}_{i=1}^N$ be a
collection of functions $v_i(t):\IZ_+^1\to X$ with the parameter
$t\in\IZ_+^1$ and let the function $u(t):\IZ_+^1\to X$ be defined
by the following (non-autonomous) relation: %
\beq{e:local-action}{
  u(t+1) := \Phi(\tau_0 u(t), \tau_1 v_1(t), \tau_2 v_2(t), \dots) ,}%
where $u(0):=u_0\in X$, the maps $\tau_i:X\to X$ satisfy the
inequality $\Lambda_{\tau_i}\le\Lambda_{\map}$ for all $i$, and
the map $\Phi:X^{N+1}\to X$ satisfy the inequality %
\beq{e:loc-Lip}{\rho(\Phi(z_0,z_1,\dots), ~\Phi(z'_0,z'_1,\dots))
               \le \Lambda_\Phi \sum_{i=0}^{N}\rho(z_i,z_i')}
for any two sequences $\{z_i\}, ~\{z'_i\}$. Then whenever
$\Lambda_\Phi \Lambda_\map<1$ there exists a function $\hat
u(t):\IZ_+^1\to X$ (depending on $v_i$) such that for any $u_0\in
X$ the solution $u(t)$ of the equation (\ref{e:local-action}) with
$u(0)=u_0$ converges with time at exponential rate to $\hat
u(t)$, i.e. $\rho(u(t), \hat u(t))\le C\gamma^t$ for all $t\ge0$
and some $\gamma<1, ~ C<\infty$. Moreover, if there exists a
vector $\vec{v}\in X^{N}$ such that %
$\max_{1\le i\le N}\rho(v_i(t),v_i) \le \hat C \hat\gamma^t$ for
all $t\ge0$ and some $\hat\gamma<1, ~ \hat C<\infty$, then the
limit $\hat u:=\lim_{t\to\infty}\hat u(t)$ exists as well.
\end{lemma}
\proof Choose two different initial values $u_0, u_0'\in X$ and
denote the corresponding solutions of the equation
(\ref{e:local-action}) by $u(t)$ and $u'(t)$ respectively. Then %
\bea{ \rho(\!\!\!\!\a u(t+1), u'(t+1)) \\
 \a= \rho(\Phi(\tau_0 u(t), \tau_1 v_1(t), \tau_2 v_2(t), \dots),
  ~\Phi(\tau_0 u'(t), \tau_1 v_1(t), \tau_2 v_2(t), \dots))  \\
\a\le \Lambda_\Phi \rho(\tau_0 u(t), \tau_0 u'(t))
\le \Lambda_\Phi \Lambda_\map \rho(u(t), u'(t)) \le \dots
\le (\Lambda_\Phi \Lambda_\map)^t \rho(u_0, u'_0)
\toas{t\to\infty}0 ,}%
which proves the first assertion with %
$C:=\max_{x,y\in X}\rho(x,y)$ and $\gamma:=\Lambda_\Phi
\Lambda_\map<1$.

To prove the second assertion consider %
\bea{ \rho(u(t+1), u(t))
 \a= \rho(\Phi(\tau_0 u(t), \tau_1 v_1(t), \tau_2 v_2(t), \dots),
       ~\Phi(\tau_0 u(t-1), \tau_1 v_1(t-1), \tau_2 v_2(t-1), \dots)) \\
\a\le  \Lambda_\Phi \Lambda_\map
     \left(\rho(u(t),u(t-1))
         + \sum_{i=1}^{N}\rho(v_i(t),v_i(t-1))\right) .}%
Using the notation $\gamma:=\Lambda_\Phi \Lambda_\map<1$ and
continuing this estimate backward in time we get %
\bea{ \rho(u(t+1), u(t))
   \a\le \gamma^2 \rho(u(t-1),u(t-2))
     + \gamma^2 \sum_{i=1}^{N}\rho(v_i(t-1),v_i(t-2)) \\
     \a + \gamma \sum_{i=1}^{N}\rho(v_i(t),v_i(t-1)) \\
   \a\le \dots
   \le \gamma^t \rho(u(1),u(0))
     + \sum_{k=1}^{t}\gamma^k \sum_{i=1}^{N}\rho(v_i(t-k+1),v_i(t-k)) \\
   \a\le \gamma^t \rho(u(1),u(0))
     + 2N\hat{C}\sum_{k=1}^{t}\gamma^k \hat\gamma^{t-k+1} \\
   \a\le \gamma^t \rho(u(1),u(0))
     + 2N\hat{C}t~(\max(\gamma,\hat\gamma))^{t+1} 
     \toas{t\to\infty}0 }%
at exponential rate (since $\max(\gamma,\hat\gamma)<1$), 
which yields the result. \qed

\n{\bf Proof} of Theorem~\ref{t:LRA}. Choose any element
$i\in\cB$. According to the definition of the map $\cmap_{|\cB}$
and to the property of the unidirectional interaction the value
$(\cmap_{|\cB}\vec{x})_i$ depends only on $x_i$ and the values at
sites with smaller indices according to the ordering given by the
function $\cN$.

Our strategy is to apply Lemma~\ref{l:contraction} consecutively
(according to the enumeration $\cN$) to all elements of the
box $\cB$ starting from the `smallest' one (say $i\in\cB$) having
the number 1. In that case the corresponding value
$(\cmap_{|\cB}\vec{x})_i$ depends only on $x_i$ and the values at
sites outside of the box $\cB$ given by boundary conditions.
Thus $x_i(t)$ is defined by the relation of type (\ref{e:local-action})
with $\Lambda_\Phi:=\Lambda_\cI$. Hence by Lemma~\ref{l:contraction}
$x_i(t)$ converges to a certain limit $\hat{x}_i$ as
time $t\to\infty$.

Then we apply Lemma~\ref{l:contraction} to the site corresponding
to the second smallest value of $\cN$, etc. Observe that the
exponential rate of convergence needed to apply the second
assertion of Lemma~\ref{l:contraction} follows from the first one.

It is important that in any given step of the procedure the values
at sites with smaller indices (according to the enumeration $\cN$)
whose entries enter into the formula (\ref{e:local-action}) and
belong to the corresponding collection $\Gamma$ do not depend on
the value of the site under consideration.
\qed

\begin{corollary} Under the assumptions of Theorem~\ref{t:LRA}
the limit solution $\vec{\hat{x}}$ depends continuously on the
boundary conditions.
\end{corollary}

Let us study how a limit solution depends on changes in
boundary conditions. To simplify calculations we consider here
only the case of linear one-dimensional interactions given
by (\ref{e:lin-unidir}) and identical one-dimensional maps
$\map_i\equiv\map$ acting on a unit circle. We assume also
that the boundary conditions are homogeneous, i.e.
$x_i=v ~ \forall i\in\cB^c$, and in some neighborhood $U$ of
the point $v\in X$ the map $\map$ is linear and expanding, i.e.
$\map x=ax+b$ with $|a|>1$ for all $x\in U$. Following the same
argument as in the proof of Theorem~\ref{t:LRA} we start with 
the site corresponding to the smallest positive value of the 
function $\cN$. Denote the value of the corresponding limit 
solution (existing by Theorem~\ref{t:LRA}) at this site by $u$. 
Then $u$ satisfies the relation
$$ u = (1-c) \map u + c\map v ~.$$
Substituting the local linear representation of $\map$ we get
$$ u = (1-c)(au+b) + c(av+b) ~,$$
or
$$ u = \frac{cav+b}{ca-(a-1)} ~.$$
Choose another homogeneous boundary condition $v'\in U$. Then
the corresponding limit solution can be written as
$$ u' = \frac{cav'+b}{ca-(a-1)} ~.$$
Therefore
$$ |u-u'| = \frac{ca}{ca-(a-1)}\cdot|v-v'| $$
and the constant $\frac{ca}{ca-(a-1)}>1$ since $|a|>1$ and 
$(1-c)a\le\La_\cI\La_\map<1$. 
Note that we need to assume here that $u,u'\in U$, so the 
neighborhood $U$ should be large enough.

Continuing this procedure for sites with larger $\cN$-indices
we see that the distance between limit solutions corresponding
to the homogeneous boundary conditions $v$ and $v'$ grows
exponentially with distance from the boundary of the box
(at least until one of limit solutions leaves $U$).

\section{Proof of Theorem~\ref{t:equiv}}\label{s:enum}

First assume that the interaction in $\cB$ is unidirectional
but the graph $G(\cB)$ has a cycle contrary to our claim.
Then when we are moving along this directed cycle the enumeration
function $\cN(\cdot)$ can only grow which contradicts
to the existence of the cycle.

Assume now that the directed graph $G(\cB)$ has no cycles.
We construct the enumeration function $\cN(\cdot)$ as follows.
We shall say that a vertex $i\in\cB$ is a {\em starting} one
if there is no $j\in\cB$ such that the edge $(j,i)\in G$, and
{\em free} if there is no $j\in\cB$ such that $(i,j)\in G$.

Clearly due to the absence of cycles in $G$ we deduce that the set
of starting vertices is non empty.

We start the enumeration procedure by choosing some starting
vertex say $i\in\cB$ and setting $\cN(i):=1$.

Assume now that we already have enumerated $M$ vertices
from $\cB$. Then we normalize the function $\cN$ in such a way
that for each $k=1,2,\dots,M$ its $k$-th largest value becomes
equal to $k$ (i.e. it is an exact numbering).

Choose an arbitrary path $i_1,i_2,\dots,i_n$ in the graph $G$
satisfying the properties that the function $\cN$ is not defined
on all vertices of this path, the vertex $i_1$ is starting or
$\exists j\in\cB: ~ (j,i_1)\in G$ and $\cN(j)$ is well defined,
and either (a) $i_n$ is free or (b) $\exists j'\in\cB: ~
(i_n,j')\in G$ and $\cN(j')$ is well defined.

If the condition (a) is satisfied we set
$$ \cN(i_k) := M + k, ~~ k=1,2,\dots,n ,$$

otherwise if the condition (b) is satisfied we set
$$ \cN(i_k) := \cN(j')-1 + \frac{k}{n+1}, ~~ k=1,2,\dots,n .$$

After this we normalize the function $\cN$, choose a new path,
etc., until no non-enumerated elements left in $\cB$.

Our aim was to define the function $\cN$ on this path in such a way
that the unidirectional interaction property should hold with respect
to previously numbered vertices. It is straightforward to check that
the enumeration obtained on each step of our construction does
satisfy this property. \qed

\section{Discussion of boundary conditions of different types}
\label{s:bc}

The result of Theorem~\ref{t:LRA} may suggest that at least in
the case of homogeneous boundary conditions the mechanism leading
to the LRA phenomenon is a simple propagation of the boundary
conditions through the boundary into a box $\cB$. In general this
is not the case and the limit solution $\vec{\hat{x}}$ can be
rather inhomogeneous in space. However this propagation phenomenon
still may take place as we shall demonstrate in the following
result. Assume that all local maps $\map_i$ are identical and
denote this common local map by $\map$.

\begin{lemma}\label{l:propagation} Let the map $\map$ have a fixed
point $\xi\in X$ (i.e. $\map\xi=\xi$) and suppose that the
boundary conditions are fixed at the value $\xi$, i.e.
$x_i(0):=\xi ~ \forall i\in\cB^c$. Then under the conditions of
Theorem~\ref{t:LRA} for any initial data $\vec{x}_{|\cB}(0)$ the
solution $\vec{x}(t)$ will converge with time to the constant
solution with $(\vec{\hat{x}})_i\equiv\xi$ for all $i\in\cB$.
\end{lemma}
\proof First, according to Theorem~\ref{t:LRA} there exists the
only one limit solution $\vec{\hat{x}}$. On the other hand, by
definition the interaction $\cI$ preserves the ``diagonal'' and
hence any sequence with identical coordinates is preserved by $\cI$.
Therefore the vector with identical coordinates
$(\vec{\hat{x}})_i\equiv\xi$ for all $i\in\cB$ is invariant for
both the map and the interaction and is, indeed, the only solution
to our problem. \qed

It turns out that the solutions depend rather sensitively
on the way how one defines the corresponding boundary conditions.

Let us start with a very mild change in the definition of the
interaction between the sites inside of the box $\cB$ and the outer
ones. Previously we used the convention that in order to calculate
the value $(\cmap_{|\cB}\vec{x})_i(t+1)$ one first applies the map
$\vec{\map}$ at all sites and then applies the interaction $\cI$
to the result of the previous operation. Contrary to this, assume
that we apply the map only to the sites from the box $\cB$ (but
not to the outer sites) and denote the map describing the corresponding
system by $\cmap^{(0)}_{|\cB}$. From the point of view
of Theorem~\ref{t:LRA} this change does not matter much and its
claim remains valid for the modified map $\cmap^{(0)}_{|\cB}$.
Indeed, it is straightforward to check that the modified map
$\cmap^{(0)}_{|\cB}$ with certain boundary conditions
$\vec{x'}_{\cB^c}$ is equivalent to the original map $\cmap_{|\cB}$
with the boundary conditions $x_i\in\map_i^{-1}x'$ for all
$i\in\cB^c$. Of course, we assume here that the preimages do exist, 
i.e. $\map_i^{-1}x'\ne\emptyset$. 
On the other hand, the result about the propagation of homogeneous
boundary conditions (Lemma~\ref{l:propagation}) is no longer valid.
The problem is that the homogeneous vector $\vec{x}$ immediately
becomes inhomogeneous under the action of $\cmap^{(0)}_{|\cB}$.

The boundary conditions considered so far can be called {\em
frozen} because they do not change in time. We shall consider now
two other types of boundary conditions. The first of them,
called {\em free}, corresponds to the situation when after each
application of the map the values at the sites outside of the box 
$\cB$ do change under $\map$, i.e. %
\beq{e:map-in-box-free}{
  \left(\cmap_{|\cB}^{({\rm free})}~\vec{x}\right)_i :=
   \function{
     \left(\left(\cI\circ\vec\map\right)\vec{x}\right)_i
           &\mbox{if } i \in \cB\\
     \map_i x_i    &\mbox{otherwise}.} }%

\begin{theorem}\label{t:LRA-modified} Let the interaction in a
finite box $\cB$ be unidirectional, the boundary conditions are
free, for any $i\in\cB$ there be a path in $G(\cB,\cI)$ from a
site outside of the box $\cB$ to $i$, and let
$\Lambda_\cI\Lambda_\map<1$. Then for any (initial) boundary
conditions $\vec{x}_{\cB^c}(0)$ there exists a vector-function
$\vec{\hat x}(t)$ such that
$$ \vec\rho(\left(\cmap_{|\cB}^{({\rm free})}\right)^t\vec{x}(0),
            ~~\vec{\hat x}(t)) \toas{t\to\infty}0 .$$
for any initial data $\vec{x}_{\cB}(0)$.
\end{theorem}

\proof Observe first that Lemma~\ref{l:contraction} has been
proven for time dependent entries. Now Theorem~\ref{t:LRA-modified}
follows essentially from the same argument as the one used
in the proof of Theorem~\ref{t:LRA}. This proves the existence
of the vector-function $\vec{\hat{x}}(t)$ which is the limit
solution (depending on time now) to the problem with free boundary
conditions. In fact the only difference to the frozen case is
that the limit solution is no longer a constant.
\qed

A completely different situation takes place if one considers
periodic boundary conditions. Of course, in this case a box $\cB$
must have a rectangular shape\footnote{or a more general polygon 
    which one can fill in by reflections the entire plane.} %
(otherwise the periodic continuation
is not well defined). Then the dynamics inside the box is completely
determined by the initial data inside the box $\cB$. However, in
distinction to the previous cases a limit solution might not be
unique. Assume for simplicity that $\map_i\equiv\map$ for all $i$
and the map $\map$ has two different periodic points $\xi$ and
$\eta$ with non-intersecting trajectories. Consider two
different initial data:
$\vec{x}_i(0):=\xi, ~ \vec{x'}_i(0):=\eta$ for $i\in\cB$.
Then the functions $\vec{x}(t):={\vec{\map}}^t\xi$ and
$\vec{x'}(t):={\vec\map}^t\eta$ describes two different
time-periodic solutions to the problem having no common limit points.
In fact, we even do not need any specific property of the interaction
here and use only the invariance of the diagonal.

\section{Stochastic perturbations of local dynamics}
\label{s:random}

Conditions under which the existence of the LRA has been proven
in Theorem~\ref{t:LRA} look very robust and certainly remain
valid under small deterministic perturbations of local
dynamics.\footnote{Namely, the vector $\vec{\hat x}$
   depends continuously on the local maps $\map_i$.} %
Our aim in this section is to show that the situation
changes drastically if instead of deterministic
perturbations one considers random ones. For a general
introduction to random perturbations of chaotic dynamical
systems we refer the reader to \cite{Bl}. However, since
the interaction map $\cI$ strictly speaking is not only
non-contracting but non-hyperbolic as well, one cannot use
general results here and everything has to be done ``by hands''.

From now on we assume that all the assumptions of
Theorem~\ref{t:LRA} are satisfied and hence the LRA
phenomenon takes place for the unperturbed system.
To simplify notation and to emphasize that the stochastic
instability of the LRA is not caused by some artificial
assumptions, throughout this section we assume
that the local maps coincide (i.e. $\map_i\equiv\map$)
and the local phase space $X$ is a unit $d$-dimensional
flat torus with the Euclidean metric on it. Note however
that all results that we shall discuss remain valid
when the maps $\map_i$ are different but satisfy
to assumptions we shall make for $\map$, and one can replace
the unit torus by a smooth compact $d$-dimensional manifold.

\OPR{A {\em stochastic perturbation} of a dynamical
system $(\map,X)$ is a Markov chain with the same phase
space $X$ defined by the {\rm transition probabilities} $Q(x,A)$
to jump from a point $x\in X$ to a Borel set $A\subset X$.
By a {\em stochastically perturbed} system we mean a
superposition of the original system and a perturbation, i.e.
a new Markov chain defined by the transition probabilities
$Q(\map x,A)$.}

\OPR{The {\em amplitude} of the stochastic
perturbation is defined as the smallest number $\ep$ such that
$Q(\map x,A)=0$ whenever $\min_{y\in A}\rho(x,y)>\ep$ for all
$x\in X$. Therefore we shall use the notation $Q_\ep$ to
emphasize that the amplitude of the perturbation is $\ep$.
We shall assume also that the transition probabilities
are translationally invariant, i.e. for any $x,z\in X$ and
a Borel set $A\subseteq X$ %
\beq{e:transl-inv}{Q(x,A)\equiv Q(x+z,A+z) ~.}  }%
Note that the phase space $X$ is a torus here and hence the 
addition is well defined.

To describe both deterministic and stochastic counterparts of
the dynamics in the same terms it is useful to consider their
action on the space of probabilistic measures $\cM(X)$ on $X$.
To distinguish between the pointwise action and the action on
measures we shall use the upper index star for the latter
(e.g. $\*\map$ stands for the action on measures for the map $\map$).
These objects are defined as follows:
$$ \*\map\mu(A) := \mu(\map^{-1}A), \qquad
   \*Q_\ep\mu(A) := \int Q_\ep(x,A) d\mu(x) $$
for any probabilistic measure $\mu\in\cM(X)$ and any Borel
set $A\subseteq X$. According to the previous definition an
action of the stochastically perturbed system on measures
is a product of these operators: $\*Q_\ep\*\map$,
while the action of the stochastically perturbed CML on
measures can be written as
$$ \*\cmap_\ep := \*\cI \*{\vec{Q}}_\ep \*{\vec\map}
               ~:~~ \cM(\vec{X}) \to \cM(\vec{X}) ~,$$
where $\*{\vec{Q}}_\ep$ stands for the direct product of
the operators $\*Q_\ep$ corresponding to local stochastic
perturbations.

\OPR{A solution $\vec\mu_\ep\in\cM(\vec{X})$ to the equation
$\cmap_\ep^* \vec\mu = \vec\mu$
is called an {\em invariant measure} for the stochastically
perturbed system.}

In distinction to our previous results which were concerned
mainly with the stability type analysis we will need a number of
assumptions to prove the instability against stochastic
perturbations. Let us introduce the following notation.
For $\mu\in\cM(X)$ denote by $\supp(\mu)$ a {\em support}
of the measure $\mu$ which defined as the intersection of all
closed subsets of $X$ of full $\mu$-measure. We define also the
collection of {\em projection operators} $\pi_k$ acting
on the space $\cM(\vec{X})$ where the index $k\in\IZ^d$
corresponds to the $k$-th component of the direct product space
$\vec{X}$, i.e. $\pi_k\vec\mu$ stands for the projection of the
measure $\vec\mu\in\cM(\vec{X})$ to the $k$-th component of the
direct product of $\cM(X)$.

We shall use the notation $B_{r}(x)$ for the ball of radius $r$
centered at a point $x\in\cB(X)$ and $\delta_{x}$ for the
$\delta$-measure at point $x$. To control a spread out of an
initial measure $\vec\mu\in\cM(\vec{X})$ under the action of the
dynamics we introduce the following notion and a list of
assumptions.

\OPR{By {\em upper/lower spreads} $S_{\pm}(\mu)$ of a measure
$\mu\in\cM(X)$ we shall mean the radius of the smallest ball
which contains support of the measure $\mu$ in the case of
$S_{+}(\mu)$, and the radius of the largest ball contained in the
support of the measure $\mu$ in the case of $S_{-}(\mu)$.}

The first of those functionals characterizes the maximal distance
between points in the support of a measure, while the second one
controls its concentration properties.

We assume that there are constants
$\La_\map\ge\la_\map>1$ and $\la_Q,a,b,\sigma>0$ such that%
\beq{e:A1}{\la_\map\rho(x,y) \le \rho(\map x ,\map y)
                             \le \La_\map \rho(x,y)
           \qquad \forall x,y: ~ \rho(x,y)\le\sigma ~,}%
\beq{e:A2}{\supp(\*Q_{\ep}\delta_{x})\supseteq B_{\la_Q\ep}(x)
           \qquad \forall x\in X ~,}%
\beq{e:A3}{S_{+}(\pi_{\ell}\*\cI\vec\mu)
             \ge a S_{+}(\pi_{\ell}\vec\mu)
               + b\sum_{\ell'}S_{+}(\pi_{\ell'}\vec\mu) }%
for any $\ell\in\cB$. The summation $\sum_{\ell'}$ here and in
the sequel is taken over all sites $\ell'\in\Gamma(\ell)$ such
that $(\ell',\ell)\in G$. We assume also that the maximal degree
of vertices of the graph $G$ is bounded by $K<\infty$.

\begin{theorem}\label{t:stoch-inst}
Let the local maps, stochastic perturbations and the interactions
satisfy the assumptions (\ref{e:A1}), (\ref{e:A2}), (\ref{e:A3})
and we consider a finite box $\cB$ with certain frozen boundary
conditions. Assume additionally that %
\beq{e:expansion}{a\la_\map<1, ~(a+b)\la_\map\ge1 .}%
Then the stochastically perturbed system in the box $\cB$
admits an absolutely continuous invariant measure
$\vec\mu_\ep\in\cM(\vec{X})$, and there exists a constant
$\gamma>0$ (which do neither depend on $\ep>0$ nor on $\cB$)
such that
$$ S_-(\pi_\ell\vec\mu_\ep)\ge\gamma L(\ell)\ep $$
for each site $\ell\in\cB$ and small enough $\ep>0$.
\end{theorem}

In other words the lower spread of the invariant measure grows at
least linearly with the `distance' from the boundary $L(\ell)$. 
Note the difference to a `diffusive' growth which one might 
expect from random perturbations. 

Let now demonstrate that the assumptions 
(\ref{e:A1} -- \ref{e:expansion}) are natural. 
First note that the assumption (\ref{e:A1})
together with the assumption $\la_\map>1$ imply only that the map
$\map$ is locally expanding. Thus, e.g. the dyadic map $\map
x:=2x~{\rm mod}1$ defined on a unit circle $X$ considered as the
interval $[0,1)$ with the metric $\rho$ on it defined as the
length of the shortest arc between points will satisfy them if we
set $\la_{\map}=\La_\map=2$ and $\sigma=1/3$. On the other hand,
the assumption $a\la_\map<1$ is necessary to be consistent with 
the LRA phenomenon in an unperturbed system (Theorem~\ref{t:LRA}).

To understand the reason behind the assumption~(\ref{e:A3}) consider
the simplest linear one-dimensional interaction defined by the
relation (\ref{e:lin-unidir}), for which we have $a=1-c$ and
$b=c$. Then (\ref{e:A3}) and the inequality $(a+b)\la_{\map}\ge1$ 
clearly hold true.

The assumption (\ref{e:A2}) gives the information about the local
stochastic perturbations of a point mass and clearly is satisfied
in the case of independent local stochastic perturbations
uniformly distributed in an $\ep$-neighborhood of a point with
$\la_{Q}=1$.

\bigskip

\n{\bf Proof of Theorem~\ref{t:stoch-inst}}. Since the
construction we shall apply is somewhat involved we start with
the main scheme of the proof. Similar to the construction of the
limit solution in the proof of Theorem~\ref{t:LRA} we shall
construct the invariant measure $\vec\mu_\ep$ step by step
starting from the first site according to the enumeration $\cN$
and adding each time only one site according to this
enumeration.\footnote{In fact, the limiting
   solution constructed in the proof of Theorem~\ref{t:LRA}
   describes the support of the invariant measure $\vec\mu_{0}$
   of the unperturbed system.} %
Using that the dynamics on already considered sites do not depend
on the dynamics of the new one (but not vice versa) we may (and
will) consider the interaction with those sites as an additional
stochastic perturbation. On this stage our main goal is to show
that the action on measures of the stochastically perturbed
system is a contraction in a suitable metric in the space of
measures. To this end we shall show that the expansion rate of
the operator $\*\map$ is proportional to $\La_\map$, while for
stochastic perturbation $\*Q$ it is close to $1$. The key point
is to show that the expansion rate of projections of the operator
$\*\cI$ are proportional to $\La_\cI$, which together with the
assumption $\La_\map\La_\cI<1$ guarantees the contraction of the
operator $\*\cmap_\ep$ and hence the convergence in time to the
unique invariant measure. On the next step of the construction
we shall obtain an upper bound for the upper spread of the
projection of the invariant measure $\vec\mu_\ep$ to the $\ell$-th
subsystem as a function of the `distance' to the boundary $L(\ell)$
and then using it we shall get the lower bound for the lower
spread of this projection.

We introduce the so called Hutchinson metrics in the space
of probabilistic measures $\cM(X)$ defined by the relation
$$ \*\rho(\mu,\nu)
  := \sup_{\La_\phi\le1}(\int\phi d\mu - \int\phi d\nu) ~,$$
where the supremum is taken over all continuous functions
$\phi:X\to\IR^1$ satisfying the assumption that their Lipschitz
constants $\La_\phi$ (defined by the relation
(\ref{e:contr-const})) do not exceed 1.

\begin{lemma} Under the assumptions of Theorem~\ref{t:stoch-inst}
for any pair of measures $\mu,\eta\in\cM(X)$ we have %
\beq{e:exp-*-T}{\*\rho(\*\map\mu,\*\map\nu)
                  \le\La_\map^\alpha \*\rho(\mu,\nu) ~,}
\beq{e:exp-*-Q}{\*\rho(\*Q_\ep\mu,\*Q_\ep\nu)
                \le \*\rho(\mu,\nu) ~,} %
\beq{e:exp-*-I}{\*\rho(\pi_\ell\*\cI\vec\mu, \pi_\ell\*\cI\vec\nu)
                \le \La_\cI\*\rho(\mu,\nu) ~,} %
where $\vec\mu, \vec\nu\in\cM(\vec{X})$ are arbitrary measures
satisfying the assumption that $\pi_\ell\vec\mu=\mu$ and
$\pi_\ell\vec\nu=\nu$.
\end{lemma}
\proof By definition
\bea{ \*\rho(\*\map\mu,\*\map\nu)
 \a= \sup_{\La_\phi\le1}
        (\int\phi d\*\map\mu - \int\phi d\*\map\nu) \\
 \a= \sup_{\La_\phi\le1}
        (\int\phi\circ\map d\mu - \int\phi\circ\map d\nu) ~.}%
Since $\rho(\map x,\map y)\le\La_\map\rho(x,y)$ for any pair of
points $x,y\in X$ we have
$$ |\phi\circ\map(x) - \phi\circ\map(y)|
 \le \La_\phi \rho(\map x,\map y)
 \le \La_\phi \La_\map \rho(x,y) ~.$$
Thus setting $\psi:=\La_\map^{-1}\phi\circ\map$ we get
$\La_\psi\le\La_\phi$ and hence
$$ \*\rho(\*\map\mu,\*\map\nu)
   \le \La_\map\sup_{\La_\psi\le1}
       (\int\psi d\mu - \int\psi d\nu)
     = \La_\map \*\rho(\mu,\nu) ~,$$
which proves inequality~(\ref{e:exp-*-T}).

Let us deal now with stochastic perturbations. %
\bea{ \*\rho(\*Q_\ep\mu,\*Q_\ep\nu)
 \a= \sup_{\La\phi\le1}
        (\int \phi d\*Q_\ep\mu - \int \phi d\*Q_\ep\nu) \\
 \a= \sup_{\La\phi\le1}
        (\int Q_\ep\phi d\mu - \int Q_\ep\phi d\nu) ~.}%
Using the translation invariance of $Q_\ep(\cdot,\cdot)$ we get %
\bea{ |Q_\ep\phi(x) - Q_\ep\phi(y)|
 \a= |\int \phi(z) Q_\ep(x,dz) - \int \phi(z) Q_\ep(y,dz)| \\
 \a= |\int \phi(z) Q_\ep(x,dz) - \int \phi(z-x+y) Q_\ep(y,dz-x+y)| \\
 \a= |\int \phi(z) Q_\ep(x,dz) - \int \phi(z-x+y) Q_\ep(x,dz)| \\
\a\le \int |\phi(z) - \phi(z-x+y)|\cdot Q_\ep(x,dz) %
\le \La_\phi \rho(x,y) ~.} %
Thus setting $\psi:=Q_\ep\phi$ we obtain
$$ \*\rho(\*Q_\ep\mu,\*Q_\ep\nu)
   \le \sup_{\La\psi\le1} (\int\psi d\mu - \int\psi d\nu)
     = \*\rho(\mu,\nu) ~,$$
which proves the inequality~(\ref{e:exp-*-Q}).

To prove the remaining inequality~(\ref{e:exp-*-I}) we introduce
the following notation. For a given $x\in X$ let
$\vec{x}_{(x_\ell=x)}\in\vec{X}$ mean a vector whose $\ell$-th
coordinate coincides with $x$, and for
$\mu_\ell:=\pi_\ell\vec\mu\in\cM(X)$ denote by $\mu_\ell^\bot$
the projection of the measure $\vec\mu\in\cM(\vec{X})$ to all
subsystems except the $\ell$-th one (i.e. we integrate over the
$\ell$-th coordinate). Then
$$ \int \phi(x) d\pi_\ell\*\cI\vec\mu(x)
 = \int\left(
       \int \phi((\cI \vec{x}_{(x_\ell=x)}))
            ~d\mu_\ell^\bot(\vec{x}_{(x_\ell=x)})
      \right) d\mu_\ell(x) .$$
Consider a function
$$ \psi(x) := \La_\cI^{-1}\int \phi((\cI \vec{x}_{(x_\ell=x)}))
            ~d\mu_\ell^\bot(\vec{x}_{(x_\ell=x)}) .$$
Then
$$ \La_\cI\cdot|\psi(x) - \psi(y)|
\le \int |\phi((\cI \vec{x}_{(x_\ell=x)}))
        - \phi((\cI \vec{x}_{(x_\ell=y)}))|
        ~d\mu_\ell^\bot(\vec{x}_{(x_\ell=x)}) ,$$
since
$\mu_\ell^\bot(\vec{x}_{(x_\ell=x)})
     = \mu_\ell^\bot(\vec{x}_{(x_\ell=y)})$ because the measure
$\mu_\ell^\bot$ does not depend on $x_\ell$.

On the other hand, by the inequality~\ref{e:lip-int} we have
$$ \rho(\cI \vec{x}_{(x_\ell=x)}, ~ \cI \vec{x}_{(x_\ell=y)})
\le \La_\cI\rho(x,y) ~,$$ %
which implies that
$$ |\psi(x) - \psi(y)|
\le \La_\cI^{-1}\La_\cI\rho(x,y)
    \int d\mu_\ell^\bot(\vec{x}_{(x_\ell=x)})
  = \rho(x,y) .$$
Therefore using the same argument as in the previous
estimates we obtain the inequality~(\ref{e:exp-*-I}). \qed

For a given $\ell\in\cB$ the set $\Gamma(\ell)$ consists of a
finite number of elements connected to $\ell$ (directly or
indirectly\footnote{through some other elements of $\cB$}).
Consider a subsystem consisted only of the sites belonging to
$\Gamma(\ell)$. The restricted dynamics of this subsystem does
not depend on the dynamics of the other parts of the system.

Collecting inequalities (\ref{e:exp-*-T}), (\ref{e:exp-*-Q}),
(\ref{e:exp-*-I}) we see that the operator corresponding to the
stochastically perturbed system is strictly contracting in the
$\*\rho$ metric, which (due to the compactness of this space)
immediately implies that for any initial measure
$\vec\mu\in\cM(\vec{X})$ we have:
$$ \*\rho(\pi_\ell(\*{\cmap}_\ep)^t\vec\mu, \pi_\ell\vec\mu_\ep)
   \toas{t\to\infty}0 ~.$$

Going through all sites of the box $\cB$ we construct the invariant
measure $\vec\mu_\ep$ for the entire stochastically perturbed
system restricted to the box $\cB$ with the given frozen boundary
conditions.

Now we are ready to obtain estimates from above for the upper
spread $S_{+}(\pi_{\ell}\vec\mu_\ep)$. It is straightforward to
check that under (\ref{e:A1}) for any measure $\mu\in\cM(X)$
whose upper spread satisfies the inequality
$S_{+}(\mu)\le\sigma/\La_\map$ we have %
\beq{e:A1+}{S_{+}(\map^*\mu) \le \La_\map S_{+}(\mu) ~.} %
The relation (\ref{e:lip-int}) implies that if
$S_{+}(\pi_{\ell}\mu),S_{+}(\pi_{\ell'}\mu)\le\sigma/\La_\map$ then%
\beq{e:A3+}{S_{+}(\pi_\ell \cI^*\vec\mu)
           \le \La_{\cI}S_{+}(\pi_\ell \vec\mu)
             + \La_{\cI} \sum_{\ell'}S_{+}(\pi_{\ell'} \vec\mu) .~}%
Since the amplitude of the stochastic perturbations
is bounded from above by $\ep>0$, one has %
\beq{e:A4+}{S_{+}(Q_\ep^*\mu) \le S_{+}(\mu) + \ep ~.}%
The dynamics at the site $\ell'$ does not depend on the
dynamics at the site $\ell$ (but not vice versa) therefore
we may choose an initial measure $\vec\mu\in\cM(\vec{X})$
satisfying the following conditions: %
$$ S_{+}(\pi_\ell\vec\mu)\le\sigma/(K\La_\map), \quad
   \pi_k\vec\mu\equiv\pi_k\vec\mu_\ep
   ~~ \forall k: ~ \cN(k)<\cN(\ell) .$$
Then by the inequalities (\ref{e:A1+}, \ref{e:A3+}, \ref{e:A4+})
we have: %
\beq{e:spread+}{S_{+}(\pi_\ell \cmap_\ep^*\vec\mu)
        \le \La_{\cI}(\La_\map S_{+}(\pi_\ell\vec\mu)+\ep)
          + \La_{\cI}\sum_{\ell'}
            (\La_\map S_{+}(\pi_{\ell'} \vec\mu)+\ep) ~.}
Denote $q:=\La_{\cI}\La_\map, ~ p:=\cI(K+1)$. Applying the
inequality~(\ref{e:spread+}) $t$ times and using that
$S_{+}(\pi_{\ell'}\cmap_\ep^*\vec\mu)
 \equiv S_{+}(\pi_{\ell'}\vec\mu)$ we get
$$ S_{+}(\pi_\ell (\cmap_\ep^*)^t\vec\mu)
   \le q^t S_{+}(\pi_\ell\vec\mu)
     + (q \sum_{\ell'}S_{+}(\pi_{\ell'}\vec\mu) + p\ep)
       \cdot \frac{1-q^t}{1-q} .$$
Now passing to the limit as $t\to\infty$ and using that
$q=\La_{\cI}\La_\map<1$ (by the condition of Theorem~\ref{t:LRA})
we obtain the desired estimate from above:
$$ S_{+}(\pi_\ell\vec\mu_\ep)
   \le \frac{q}{1-q} \sum_{\ell'}S_{+}(\pi_{\ell'}\vec\mu_\ep)
     + \frac{p\ep}{1-q} ~.$$
Thus the upper spread of the $\ell$-th projection of the
invariant measure cannot exceed a constant in power $L(l)$
multiplied by $\ep$. Observe that for each $t\ge1$ the inequality
$D(\pi_\ell(\cmap_\ell^*)^t\vec\mu)\le\sigma$ holds automatically
since $q<1$.

Now whence we have established the control over the upper spread
of $\vec\mu_\ep$ we will obtain an estimate from below for the
lower spread of this measure.

Under (\ref{e:A1}) for any measure $\mu\in\cM(X)$ with
$S_-(\mu)\le\sigma/\La_{\map}$ we have: %
\beq{e:A1'}{S_-(\*\map\mu) \ge \la_{\map}S_-(\mu) ~.} %
The relation~(\ref{e:A2}) implies %
\beq{e:A2'}{S_-(\*Q_{\ep}\mu) \ge S_-(\mu) + \la_{Q}\ep ~.} %
Making use of the inequalities (\ref{e:A1'}), (\ref{e:A2'}) and
(\ref{e:A3}) for any measure $\mu\in\cM(X)$ with
$S_-(\mu)\le\sigma/\La_{\map}$ and a measure
$\vec\mu\in\cM(\vec{X})$ with $\pi_{\ell}\vec\mu=\mu$ we get: %
\beq{e:A4}{ S_-(\pi_\ell \cmap_\ep^*\vec\mu)
 \ge a(\la_\map S_-(\pi_\ell\vec\mu) +\la_Q\ep)
   + b(\la_\map \sum_{\ell'}S_-(\pi_{\ell'}\vec\mu) + \la_Q\ep) ~.}%
Denote $\alpha:=a\la_\map<1$, $\beta:=b\la_\map \ge 1-\alpha$ (by
the assumption~(\ref{e:expansion})) and $\kappa:=(a+b)\la_Q$.
Now, applying~(\ref{e:A4}) to $(\cmap_\ep^*)^t$ we have
$$ S_-(\pi_\ell(\cmap_\ep^*)^t\vec\mu)
 \ge \alpha^t S_-(\pi_\ell\vec\mu)
   + (\beta \sum_{\ell'}S_-(\pi_{\ell'}\vec\mu) + \kappa\ep)
     \cdot\frac{1-\alpha^t}{1-\alpha} .$$
Passing to the limit and using that $\beta/(1-\alpha)\ge1$ we
obtain the final estimate
$$ S_-(\pi_\ell \vec\mu_\ep)
   \ge \sum_{\ell'}S_-(\pi_{\ell'} \vec\mu_\ep)
     + \frac{\kappa\ep}{1-\alpha} ~.$$
Theorem~\ref{t:stoch-inst} follows by setting
$\gamma:=\kappa\ep/(1-\alpha)$. \qed


\section{Conclusion}\label{s:conclusion}

We have shown that the absence of loops in directed connectivity
graphs of networks of dynamical elements results in appearance of
a long range action provided that interactions between local
dynamical systems (elements of a network) are sufficiently strong. 
If local systems have a regular (contractive) dynamics than the LRA
appears already for arbitrarily weak interactions (conditions of
Theorem~\ref{t:LRA} clearly hold in this case).

One expects that this restriction to the network topology is
by no means a necessary one. In principle, we see two ways
to overcome this problem. The first of them is to use some
rather general ergodic-theoretical results (see e.g. \cite{Si})
to recode interaction potentials with one-sided dependence to
the ones with the two-sided dependence. It is known that if
those potentials are homological to each other than they
generate the same Gibbs measures (corresponding to our
invariant measures). The problem here is that it is hard to
describe the set of constructed potentials corresponding to
general non-unidirectional interactions.

Another possibility is to consider non-unidirectional interactions
as small perturbations to the unidirectional ones. This approach
seems to be more promising but one inevitably needs to make a
number of additional assumptions about the unperturbed system
to ensure that the unidirectional interactions are stronger
than the perturbations.

There is yet another possibility related to the assumption that
the interaction in the network is close to a mean-field model,
i.e. a large part of the network is interacting with each its
element. The main disadvantage of this assumption is that it is
not completely clear how to extend this property for infinite
networks. 

It is worth mentioning that instead of the lattice $\IZ^d$ one 
can consider a countable ordered collection of vertices. 
The short range interactions in this context would mean that 
subsystems corresponding to vertices interact only if the 
distance between them according to the above mentioned order 
does not exceed a certain threshold. Clearly all our results 
remain valid in this setting.

The problems similar to the ones discussed in this paper were
studied in the context of the synchronization theory. More
specifically, we are talking here about a complete synchronization
of dynamical systems \cite{PRK}. The complete synchronization
means that asymptotically, as time goes to infinity, the
differences in the behavior of all elements in the network
disappear. In particular, in \cite{Wu} it has been shown that in
the case of the network described by the ordinary differential
equations whose directed connectivity graph is a tree the
synchronization phenomenon takes place under strong enough
interactions, while under the assumption that the connectivity
graph is ``strongly connected'' (i.e. the interaction is close to
a mean-field model) this has been proven in \cite{BBH} (see also
references in these papers and the recent monograph \cite{PRK}).
In these publications, in distinction to our approach, a
synchronization phenomenon is studied in closed finite networks
consisting of identical elements and the synchronized system may
demonstrate chaotic behavior.

\section*{Acknowledgements}

The research by M.B. was partially supported by Russian
Foundation for Fundamental Research and French Ministry of
Education grants. L.B. was partially supported by the NSF grant
DMS - 0140165 and by the Humboldt foundation.

\end{document}